\documentstyle[11pt, amscd, amsmath, latexsym, setspace, amsfonts, amssymb, verbatim]{article}

\newtheorem{th}{Theorem}[section]
\newtheorem{de}[th]{Definition}

\newtheorem{lem}[th]{Lemma}
\newtheorem{co}[th]{Corollary}
\newtheorem{re}[th]{Remark}

\newtheorem{pro}[th]{Proposition}

{

\newcommand{\espa}{\qquad}

\title{On best proximity points in metric and Banach spaces}
\author{Rafa Esp\'{\i}nola \& Aurora Fern\'andez-Le\'on\footnote{Both authors were partially supported by the Ministery of Science
            and Technology of Spain, Grant BFM 2000-0344-CO2-01 and La Junta de Antaluc\'{\i}a project
        FQM-127.}}

\begin{document}

\date{}

\maketitle





\begin{abstract} In this paper we study the existence and uniqueness of best proximity points of cyclic contractions as well as the convergence of iterates to such proximity points. We do it from two different approaches, leading each one of them to different results which complete, if not improve, other similar results in the theory. Results in this paper stand for Banach spaces, geodesic metric spaces and metric spaces. We also include an appendix on CAT$(0)$ spaces where we study the particular behavior of these spaces regarding the problems we are concerned with.

MS Classification: 54H25, 47H09
\end{abstract}

\section{Introduction}
Let $A$ and $B$ be two nonempty closed subsets of a complete metric space $X$. Consider a mapping $T:A\cup B\to A\cup B$ such that
$$
T(A)\subseteq B\hskip.5cm {\rm and}\hskip.5cm T(B)\subseteq A
$$
with the additional condition that there exists $k\in (0,1)$ such that
$$
d(Tx,Ty) \le kd(x,y) \text{ for all $x\in A$ and $y\in B$},
$$
then $A\cap B\neq \emptyset$ and $T$ has a unique fixed point in $A\cap B$.

In \cite{elve, elvek, espi, suzu} a generalization of this situation was studied under the assumption of $A\cap B=\emptyset$. More precisely, in \cite{elve, suzu} it was assumed that there exists $k\in (0,1)$ such that
\begin{equation}\label{contraction}
d(Tx,Ty)\le kd(x,y)+(1-k)\;{\rm dist}(A,B)
\end{equation}
for all $x\in A$ and $y\in B$ to obtain existence, uniqueness and convergence of iterates to the so-called best proximity points; that is, a point $x$ either in $A$ or $B$ such that $d(x,Tx)={\rm dist}(A,B)$. This was first studied in \cite{elve} for uniformly convex Banach spaces, see also \cite{kireve} for more on related topics. Then, in \cite{suzu}, the property UC (see Section 3 for definition) was introduced for a pair $(A,B)$ of subsets of a metric space so a result on existence, uniqueness and convergence of iterates stands (Theorem \ref{suzuki} in Section 2) in general metric spaces. Since, as it is also proved in \cite{suzu}, property UC happens for a large collection of pairs of subsets of uniformly convex Banach spaces, Theorem \ref{suzuki} actually contains the main theorem of \cite{elve} (Theorem 3.10 in \cite{elve}) as a particular case. Property UC was even proved, in \cite{suzu}, to happen outside the setting of uniformly convex Banach spaces. In fact, this was obtained for UCED (uniformly convex in every direction) Banach spaces and strictly convex Banach spaces but in both cases under the very strong condition (see Theorem \ref{compacidad} in Section 3) of one of the sets to be of compact closure. In this work we first introduce a new property, the so-called property WUC, which is proved to happen under far less restrictive conditions than where property UC seems to reach, and then an existence, uniqueness and convergence theorem is proved for pairs of sets verifying property WUC. Second, we focus the same problem from the new approach suggested by one of the authors in \cite{espi} to obtain still new results on the same problem. As a result, a partial answer in the positive is given to a question raised in \cite{elve}.

The work is organized as follows: in Section 2 we introduce most of the definitions, notations and previous results we will need. In Section 3 we look for weaker conditions than property UC. We introduce properties WUC and W-WUC and show that similar results to those in \cite{suzu} hold under conditions which are easier to verify. In Section 4, we approach the same problem by the introduction of a semimetric. This is applied in a successful way by showing that the mappings verifying the contractive condition (\ref{contraction}), under suitable assumptions, are contractions with respect to a certain semimetric. We finish this work with a remark on CAT$(0)$ spaces. In \cite{elvek} it was shown that when the ambient space is a Hilbert space then the kind of mappings we are dealing with actually behave as nonexpansive ones. In \cite{espi} it is shown that the semimetric there defined coincides with the metric of the ambient space when this is a Hilbert space. Our remark on CAT$(0)$ spaces, in a certain sense the nonlinear counterparts of Hilbert spaces, states that something similar happens in these spaces.

\section{Preliminaries}\label{preli}

In this section we compile the main concepts and results we will
work with along this paper. We begin with some basic definitions
and notations that are needed. Let $(X,d)$ be a metric space and let
$A$ and $B$ be two subsets of $X$. Define
\begin{align*}
{\rm dist}(x,A)= & \inf \{d(x,y) : y \in A\};\\
P_A(x)= &\{y \in A : d(x,y)={\rm dist}(x,A)\};\\
\text{dist}(A,B)= &\inf \{d(x,y) : x\in A,y\in B\};\\
\text{diam} (A)= &\sup \{ d(x,y): x,y\in D\}.
\end{align*}
Recall that the set $A$ is said to be {\it a Chebyshev set with respect to $B$} if $P_A(x)$ is a singleton
for any $x\in B$.

A metric space $(X,d)$ is said to be a
{\it geodesic space ($D$-geodesic space, respectively)} if every two points $x$ and $y$ of $X$ (with
$d(x,y)\le D$) are
joined by a geodesic, i.e, a map $c:[0,l]\subseteq {\mathbb R}\to
X$ such that $c(0)=x$, $c(l)=y$, and
$d(c(t),c(t^{\prime}))=|t-t^{\prime}|$ for all $t,t^{\prime} \in
[0,l]$. Moreover, $(X,d)$ is called {\it uniquely geodesic ($D$-uniquely geodesic)} if
there is exactly one geodesic joining $x$ and $y$ for each $x,y\in
X$ (with $d(x,y)\le D$). When the geodesic between two points is unique, its image
(called geodesic segment) is denoted by $[x,y]$. The midpoint $m$ in between two points
$x$ and $y$ in a uniquely geodesic metric space is the only point in $[x,y]$ such that
$d(x,m)=d(y,m)$. Any Banach space is a geodesic
space with usual segments as geodesic segments. Throughout this work we will just use geodesic metric space to refer to a uniquely geodesic space since all our geodesic spaces will be uniquely geodesic.

A very important class of geodesic metric spaces are the CAT$(k)$ spaces, that is, metric spaces of
curvature uniformly bounded above by $k$. These spaces have been the object of a lot of interest by many researches
and we will get back to them, especially to CAT$(0)$ spaces, at certain moments of our exposition. For a very thorough treatment on CAT$(k)$-spaces the reader can check \cite{brha}.


A subset $A$ of a geodesic metric space $(X,d)$ is said to be {\it convex}
if the geodesic joining each pair of points $x$ and $y$ of $A$ is contained in $A$.

We will need the notion of uniformly convex geodesic metric space (see also \cite[pg. 107]{gore}).

\begin{de} A geodesic metric space $(X,d)$ is said to be {\it uniformly
convex} if for any $r>0$ and any $\varepsilon \in
(0,2]$ there exists $\delta \in (0,1]$ such that for all $a, x,
y\in X$ with $d(x,a)\le r$, $d(y,a)\le r$ and $d(x,y)\geq
\varepsilon r$ it is the case that
$$%
d(m,a)\le (1-\delta)r%
$$%
where $m$ stands for a midpoint of the geodesic segment $[x,y]$.
A mapping $\delta:(0,+\infty)\times (0,2]\to (0,1]$ providing such a
$\delta=\delta (r,\varepsilon )$ for a given $r>0$ and $\varepsilon
\in (0,2]$ is called a {\it modulus of uniform convexity}. If
moreover $\delta$ decreases with $r$ (for a fixed $\varepsilon$) we
say that $\delta$ is a monotone modulus of uniform convexity of $X$.
\end{de}

The notion of monotone modulus of uniform convexity seems to have been
studied for first time in \cite{leus}. Of course, the usual modulus of convexity
of a uniformly convex Banach space is monotone in this sense. For more on geometry of Banach
spaces the reader can check \cite{aydo, goki, kisi}.

\begin{re}
If in the above definition we drop the uniformity conditions then
we find the notion of strict convexity. More precisely, if $X$ is a Banach space and such a $\delta$ exists for each $a,x$ and $y$ as above with $d(x,y)>0$, then we will say that
$X$ is a strictly convex Banach space. If the same condition is imposed on a geodesic metric space $X$, then we can find the spaces of nonpositive curvature in the sense of Busemann, see \cite{papa} for a detailed study on them.
\end{re}

Cyclic contractions and best proximity points are defined next.

\begin{de}
Let $A$ and $B$ be two nonempty subsets of a metric space
$X$. A map $T : A \cup B \rightarrow A \cup B$ is a cyclic
contraction map if it satisfies:
\begin{itemize}
\item [$(1)$] $T(A) \subseteq B$ and $T(B)\subseteq A$.
\item [$(2)$] There is some $k \in (0,1)$ such that $d(Tx,Ty) \leq
kd(x,y)+(1-k)\; {\rm dist}(A,B)$, for all $x \in A$ and $y \in B$.
\end{itemize}
\end{de}

\begin{re}
Notice that condition $(2)$ implies that $T$ is a relatively nonexpansive mapping, i.e., $T$ satisfies that
$d(Tx,Ty) \leq d(x,y)$ for all $x\in A$ and $y\in B$, which were the main object of study in \cite{elvek, espi}.
\end{re}

Next we define the notion of best proximity point.

\begin{de}
Let $A$ and $B$ be two nonempty subsets of a metric space
$X$. Let $T : A \cup B \rightarrow A \cup B$ such that $T(A)
\subseteq B$ and $T(B) \subseteq A$. A point $x \in A \cup B$ is
said to be a best proximity point for $T$ if
$d(x,Tx)={\rm dist}(A,B)$.
\end{de}

Existence, uniqueness and convergence of iterates to a best proximity point for cyclic contractions
have recently been studied in \cite{elve, suzu}. The goal of this work is to
find improvements of main results in these works. Next we state the main result
from \cite{suzu} (the definition of property UC is in Section 3).

\begin{th}\label{suzuki}
Let $(X,d)$ be a metric space and let $A$ and $B$ be nonempty subsets of $X$ such that $(A,B)$ satisfies the property $UC$. Assume that $A$ is complete. Let $T$ be a cyclic contraction on $A \cup B$. Then $T$ has a unique best proximity point $z$ in $A$ and $\{T^{2n}x\}$ converges to $z$ for every $x \in A$.
\end{th}

Main result in \cite{elve} states basically the same but with $A$ and $B$ nonempty closed and convex, $X$ a uniformly convex Banach space and no mention to property UC.

In \cite{espi} a new approach to relatively nonexpansive mappings lead to the fact that
such mappings, under suitable conditions, are actually nonexpansive with respect to an adequate
semimetric. In Section 4 we apply this new approach to cyclic contractions. Next we introduce the main notions and results on semimetric spaces that we will need.

\begin{de} Let $M$ be a nonempty set. A function $d : M \times M \rightarrow [0,\infty)$
is said to be a semimetric on $M$ if
\begin{itemize}
\item [(1)] $d(x,y)=0$ if, and only if, $x=y$. \item [(2)]
$d(x,y)=d(y,x)$ for any $x,y \in X$.
\end{itemize}

In this case, $(M,d)$ is said to be a semimetric space.
\end{de}

Contractions with respect to semimetrics are defined in a similar
way to contractions with respect to metrics.

\begin{de}\label{contra} Let $(X,d)$ be a semimetric space. A mapping $T : X \rightarrow X$
is said to be a contraction if there is a constant $k \in (0,1)$
such that for all $x,y \in X$
$$d(Tx,Ty) \leq kd(x,y).$$
\end{de}

The next definition will make easier to state some of our results.

\begin{de} \label{compatible} Let $X$ be a nonempty set. Let $d$ and $d_1$ be a metric and a semimetric on $X$ respectively.
We say that $d$ and $d_1$ are compatible on $X$ if for every
$\varepsilon
> 0$ and $x\in X$ there exist $f_x(\varepsilon)>0$ and
$g_x(\varepsilon)>0$ such that
$$%
B_d(x,f_x(\varepsilon)) \subseteq B_{d_1}(x,\varepsilon)\text{ and
}
B_{d_1}(x,g_x(\varepsilon)) \subseteq B_d(x,\varepsilon),%
$$
where $B_d(x,r)$ and $B_{d_1}(x,r)$ stand, respectively, for the
closed balls of center $x$ and radius $r$ with respect to the metric and the semimetric.
\end{de}

In \cite{jama}, different counterparts of Banach's
contraction theorem are given for semimetric spaces. We state next a particular
case of those results more adequate to our context (see Theorem 1 in \cite{jama}).

\begin{th}\label{convergencia}
Let $X$, $d$ and $d_1$ be as in the definition above with $d$ and
$d_1$ compatible. Let $T$ be a contraction on $X$ for the
semimetric $d_1$, then $T$ has a unique
fixed point $x_0$. Moreover, for any $x\in X$ the sequence $\{T^nx\}_{n=1}^{\infty}$
converges to $x_0$.
\end{th}

We finish this section introducing two geometrical properties for
Banach spaces. We begin describing property $(H)$.

\begin{de} Let $X$ be a Banach space. $X$ is said to have the property $(H)$
if for any sequence on the unit sphere of $X$, weak and norm
convergence coincide.
\end{de}

\begin{re} This property has been very extensively studied in the literature and it
is closely related to the so-called Kadec-Klee
property (KK-property, for short). For more on this topic, see \cite{aydo, goki, kisi,
PKlin}.
\end{re}

We will also need the following uniform version of the KK
property.

\begin{de}
Let $X$ be a Banach space. $X$ is said to have the property UKK (uniform Kadec-Klee property) if
for any $\varepsilon >0$ the number
$$
\eta(\varepsilon)=\inf\{ 1-\| x\|\}>0,
$$
where the infimum is taken over all points $x$ such that $x$ is a
weak limit for some sequence $\{ x_n\}$ in the unit ball of $X$
with $\| x_n-x\|\geq \varepsilon$ for all $n$.
\end{de}

Different properties of UKK Banach spaces as well as
connection among all these geometrical notions can be found in the above-mentioned
references. Let us just note here, as a matter of fact, that uniformly convex
Banach spaces are UKK spaces and so they also have property $(H)$. Both notions, uniform
convexity and property UKK, have to do with a certain rotoundity of the balls of the space.
This is obvious for uniform convexity and far less obvious for property UKK as there exist
Banach spaces which are UKK and not even strictly convex.


\section{The UC and WUC properties}

Property UC was defined in \cite{suzu} in the following way.

\begin{de} Let $A$ and $B$ be nonempty subsets of a metric space
$(X,d)$. Then $(A,B)$ is said to satisfy the property UC if for $\{
x_n\}$ and $\{x^{\prime}_n\}$ sequences in $A$ and $\{ y_n\}$ a
sequence in $B$ such that
$\lim_nd(x_n,y_n)=\lim_nd(x^{\prime}_n,y_n) = {\rm dist} (A,B)$,
then $\lim_n d(x_n,x^{\prime}_n) =0$.
\end{de}

The following proposition shows the
uniform nature of property UC.

\begin{pro}\label{equivalencia} For $A$ and $B$ nonempty subsets of a metric space $X$,
the following are equivalent:
\begin{itemize}
\item[(i)] $(A,B)$ has property UC.
\item[(ii)] For any $\varepsilon >0$ there exists $\delta >0$ such
that ${\rm diam}(A\cap B(y,{\rm dist}(A,B)+\delta))\le \varepsilon$
for any $y\in B$.
\end{itemize}
\end{pro}

{\it Proof.} First we see $(i) \Rightarrow (ii)$. Supposing the contrary implies that
there is $\varepsilon_0 >0$ such that for every $\delta=1/n$ there
exist $y_n \in B$ and $x_n, x^{\prime}_n \in A$ satisfying $d(y_n,x_n) \leq
{\rm dist}(A,B)+\frac{1}{n}$, $d(y_n,x^{\prime}_n) \leq {\rm
dist}(A,B)+\frac{1}{n}$ and $d(x_n,x^{\prime}_n) > \varepsilon_0$, which obviously contradicts property UC.

Now we prove $(ii) \Rightarrow (i)$. Let $x_n,x^{\prime}_n \in A$ and $y_n
\in B$ such that $d(y_n,x_n)$ and $d(y_n,x^{\prime}_n)$ both converge to ${\rm dist}(A,B)$ as $n\to\infty$. Then given $\varepsilon
>0$, there is $\delta >0$ such that ${\rm
diam}(A\cap B(y,{\rm dist}(A,B)+\delta))\le \varepsilon$ for any $y\in B$.
Now, it is enough to take $n_0 \in \mathbb{N}$ such that $d(x_n,y_n), d(x^{\prime}_n,y_n) \leq {\rm dist}(A,B) +
\delta$ for any $n \geq n_0$ to deduce that $d(x_n,x^{\prime}_n) \leq \varepsilon$
for $n \geq n_0$. $\Box$\vskip.5cm

In \cite{suzu} it was shown that any pair of nonempty subsets $(A,B)$ of uniformly
convex Banach spaces with $A$ convex enjoy the property UC. Next we show that something
similar can be said for uniformly convex geodesic spaces under adequate conditions
on the modulus of convexity.

\begin{pro}\label{converg1} Let $(X,d)$ be a uniformly convex geodesic metric space with
a monotone modulus of convexity $\delta(r,\varepsilon)$. Let $A$
and $B$ be two nonempty subsets of $X$ with $A$ convex. Then the
pair $(A,B)$ has property UC.
\end{pro}

{\it Proof.}
Suppose on the contrary that there exist $\{x_n\}$ and $\{x^{\prime}_n\}$ sequences in $A$,
$\{y_n\}$ in $B$ and $\varepsilon_0>0$ such that
for every $k \in \mathbb{N}$, there exist $n_k\geq k$ for which
$d(x_{n_k},x^{\prime}_{n_k}) \geq \varepsilon_0$ while $\lim_{n\to\infty} d(x_n,y_n)=
\lim_{n\to\infty} d(x^{\prime}_n,y_n)={\rm dist}(A,B)$.

There is no loss of generality in assuming that $\delta(r,\varepsilon)<1$ for $r,\varepsilon >0$
and that ${\rm dist}(A,B)>0$ since otherwise the result
follows in a trivial way. For $\gamma >
\text{dist}(A,B)$ and $\varepsilon_1=\varepsilon_0/\gamma$, choose $\varepsilon >0$  such that
$$%
\varepsilon < \min\left\{\gamma - \text{dist}(A,B),
\frac{\text{dist}(A,B)
\delta(\gamma,\varepsilon_1)}{1-\delta(\gamma,\varepsilon_1)}\right\}.
$$
Then there exists $N_0 \in \mathbb{N}$ such that if $n_k \geq N_0$,
then $d(x_{n_k},y_{n_k}) \leq \text{dist}(A,B) + \varepsilon$ and
$d(x^{\prime}_{n_k},y_{n_k}) \leq \text{dist}(A,B) + \varepsilon$.
Let $m_{n_k}$ be the mid-point of the geodesic segment
$[x_{n_k},x^{\prime}_{n_k}]$. Using the uniform convexity of $X$, we
have that
$$%
d(y_{n_k},m_{n_k}) \leq \big{(}1-\delta(\text{dist}(A,B)+ \varepsilon,
\varepsilon_1)\big{)}(\text{dist}(A,B) + \varepsilon)\leq%
$$%
$$%
\leq \big{(}1-\delta(\gamma,
\varepsilon_1)\big{)}(\text{dist}(A,B) + \varepsilon) <
\text{dist}(A,B).%
$$%

Then, for $n_k \geq N_0$
$$
d(y_{n_k},m_{n_k})< \text{dist}(A,B),
$$
which contradicts the fact that $m_{n_k}\in A$ by convexity of $A$.  $\Box$ \vskip.5cm

\begin{re} Notice that the same result remains true if the condition on the
monotonicity of the modulus of convexity is replaced by the condition of
being lower semi-continuous from the right.
\end{re}

As it was pointed in the Introduction, property UC was also shown in \cite{suzu} to happen in UCED Banach spaces and strictly convex Banach spaces but requesting $A$ is relatively compact. Regarding the assumption on the compactness of $A$ the following result from \cite{elve} is relevant.

\begin{th}
Let $A$ and $B$ be nonempty closed subsets of a metric space $(X,d)$ and let $T:A\cup B\to A\cup B$ be a cyclic contraction. If either $A$ or $B$ is boundedly compact, then there exists $x$ in $A\cup B$ with $d(x,Tx)={\rm dist}(A,B)$.
\end{th}

Notice that what we miss from Theorem \ref{suzuki} in this theorem is uniqueness and convergence of iterates. We see next that this is easy to obtain by adding the very mild condition (see the remark below to support this idea) of being $A$ a Chebyshev set for proximinal points with respect to $B$.

\begin{de}
Given $A$ and $B$ two nonempty subsets of a metric space, we say that $A$ is a Chebyshev set for proximinal points with respect to $B$ if for any $x\in B$ such that ${\rm dist}(x,A)={\rm dist}(A,B)$ we have that $P_A(x)$ is a singleton.
\end{de}

Then we can prove the following.

\begin{th}\label{compacidad}
If in the above theorem, $A$ is supposed to be boundedly compact and a Chebyshev set for proximinal points with respect to $B$, then the best proximity point $z\in A$ is unique and the sequence $\{ T^{2n}x\}$ converges to $z$ for any $x\in A$.
\end{th}

{\it Proof.} We first show it is unique. Suppose $z$ and $z^{\prime}$ are two best proximity points in $A$ with $z\neq z^{\prime}$. Then the Chebyshev condition on $A$ implies that $Tz\neq Tz^{\prime}$. Now, the relative nonexpansivity of $T$ implies that
$$
d(T^2z,Tz)\le d(z,Tz)={\rm dist}(A,B)
$$
and so, the Chebyshev condition on $A$ also implies that $z$ and $z^{\prime}$ are fixed points for $T^2$. If we write $d^*(x,y)=d(x,y)-{\rm dist}(A,B)$ then
$$%
d^*(z,Tz^{\prime})=d^*(T^2z,Tz^{\prime})%
$$%
$$%
\le k d^*(z^{\prime},Tz)=kd^*(T^2z^{\prime},Tz)\le
k^2d^*(z,Tz^{\prime}).%
$$%
Hence $d^*(z,Tz^{\prime})=0$ and so $z=z^{\prime}$.

Finally the converges of the iterates follows directly from the facts that $A$ is boundedly compact, the sequences $\{ T^{2n}x\}$ are bounded for any $x\in A$ and that $\lim d(T^{2n}x,Tz)={\rm dist}(A,B)$ for any $x\in A$. $\Box$ \vskip.5cm

\begin{re}\label{max}
Notice that the condition of being Chebyshev is a very natural one in this kind of problems. Think otherwise on the sets $A=\{ (x,0):\; x\in [0,1]\}$ and $B=\{ (x,1):\; x\in [0,1]\}$ as subsets of the plane with the maximum norm. Then any mapping $T:A\cup B\to A\cup B$ with $T(A)\subseteq B$ and $T(B)\subseteq A$ is a cyclic contraction.
\end{re}

We suggest to replace property UC with the weaker one WUC which we
define next.

\begin{de}
Let $A$ and $B$ be nonempty subsets of a metric space $(X,d)$. Then
$(A,B)$ is said to satisfy the property WUC if for any $\{
x_m\}\subseteq A$ such that for every $\varepsilon >0$ there exists
$y\in B$ satisfying that $d(x_m,y)\le {\rm dist}(A,B)+\varepsilon $
for $m\geq m_0$, then it is the case that $\{ x_m\}$ is convergent.
\end{de}

\begin{re}
Another alternative for the above definition is to ask the sequence $\{ x_m\}$ to
be Cauchy instead of convergent. It is worthwhile to note here that this is quite
a detail of a formal nature since in all our main results we always assume $A$ to be complete.
\end{re}

Next proposition gives the relation between the two mentioned
properties.

\begin{pro} Let $A$ and $B$ be nonempty subsets of a metric space
$(X,d)$ such that $A$ is complete. Suppose the pair $(A,B)$ has
property UC. Then $(A,B)$ has property WUC.
\end{pro}
{\it Proof.} Let $\{x_m\}\subseteq A$ be such that for every
$\delta >0$ there exists $y\in B$ satisfying that $d(x_m,y)\le
{\rm dist}(A,B)+\delta $ for $m\geq m_0$. It suffices to show that $\{x_m\}$ is a Cauchy sequence.
This follows directly from
 $(ii)$ of Proposition \ref{equivalencia}. $\Box$ \vskip.5cm

Next we show that property WUC implies a nonuniform version of the equivalence given
by Proposition \ref{equivalencia} for property UC. We omit its proof.

\begin{pro} Let $A$ and $B$ be nonempty subsets of a metric space $(X,d)$.
Suppose $(A,B)$ has property WUC then
$$%
\lim_{\varepsilon\to 0} {\rm diam}(A\cap B(y,{\rm
dist}(A,B)+\varepsilon))=0%
$$%
for any $y\in B$.
\end{pro}


The next propositions show that property WUC is
likely to happen in more situations than property UC. We first weaken the notion of uniform convex geodesic space with a monotone modulus of convexity.

\begin{de}
A geodesic metric space $(X,d)$ is said to be {\it pointwise uniformly convex} if for any $a\in X$, $r>0$ and $\varepsilon \in
(0,2]$ there exists $\delta=\delta (a,r, \varepsilon) \in (0,1]$ such that for all $x,
y\in X$ with $d(x,a)\le r$, $d(y,a)\le r$ and $d(x,y)\geq
\varepsilon r$ it is the case that
$$%
d(m,a)\le (1-\delta)r%
$$%
where $m$ stands for the midpoint of the geodesic segment $[x,y]$.
A mapping $\delta:X\times (0,+\infty)\times (0,2]\to (0,1]$ providing such a
$\delta$ for a given $a\in X$, $r>0$ and $\varepsilon
\in (0,2]$ is called a {\it modulus of pointwise uniformly convexity} (or just modulus of convexity when confusion cannot arise). If
moreover $\delta$ decreases with $r$ (for $\varepsilon$ and $a$) we
say that $\delta$ is a monotone modulus of pointwise uniformly convexity of $X$.
\end{de}

\begin{re}
Notice that both notions of uniform convexity coincide for Banach spaces.
\end{re}

\begin{pro} Let $(X,d)$ be a complete pointwise uniformly convex geodesic metric space with monotone modulus of
convexity. Let $A$
and $B$ be two nonempty subsets of $X$ with $A$ convex. Then the
pair $(A,B)$ has property WUC.
\end{pro}
{\it Proof.} Let $\{x_m\}$ be a sequence in $A$ such that for every $\varepsilon >0$ there exist
$y\in B$ and $m_0\in {\mathbb N}$ satisfying that $d(x_m,y)\le {\rm dist}(A,B)+\varepsilon $ for $m \geq m_0$.
Suppose $\{ x_m\}$ is not convergent, then there exists $\varepsilon_0 > 0$ such that for each $k \in {\mathbb N}$ there are $n_k,m_k \geq k$ for which $d(x_{n_k},x_{m_k}) \geq \varepsilon_0$. Thus, for $k=m_0$, we find $n_{m_0},m_{m_0} \geq m_0$ such that
$$
d(x_{n_{m_0}},y)\le {\rm dist}(A,B)+\varepsilon,\;\; d(x_{m_{m_0}},y)\le {\rm dist}(A,B)+\varepsilon
$$
and
$$
d(x_{n_{m_0}},x_{m_{m_0}}) \geq \varepsilon_0.
$$
Let $z_{m_0}$ be the mid-point in the segment $[x_{n_{m_0}}, x_{m_{m_0}}]$. Since $X$ is pointwise uniformly convex, we obtain that
$$
d(z_{m_0},y)\leq ({\rm dist}(A,B)+\varepsilon)(1-\delta),
$$
for some $\delta=\delta (y,{\rm dist}(A,B)+\varepsilon, \varepsilon_1) \in (0,1]$ as in the proof of Proposition \ref{converg1}. The contradiction follows from the fact that we can repeat this reasoning for any $\varepsilon>0$ with $y$ and $\varepsilon_0$ fixed. $\Box$\vskip.5cm

\begin{re}
This kind of modulus has been previously used for hyperbolic spaces in \cite{shafrir}.
\end{re}

\begin{pro}
Let $X$ be a UKK reflexive and strictly convex Banach space. Then,
for $A,B\subseteq X$ nonempty and convex, it is the case that
$(A,B)$ has the property WUC.
\end{pro}

{\it Proof.} Let $\{ x_n\}\subseteq A$ be as in the above proof. Suppose
$\{ x_n\}$ is not convergent. First we show that
this sequence needs to have a separated subsequence. Consider two
convergent subsequences $\{ x_{n_k}\}$ and $\{ x_{n_l}\}$ of ${x_n}$
with respective limits $x$ and $x^{\prime}$ in the closure of $A$.
For each $n\in N$ choose $y_n\in B$ such that the tales of both
subsequences are in $B(y_n,{\rm dist}(A,B)+1/n)$. Then it is clear
that
$$%
x,x^{\prime}\in \bigcap_{n\in \mathbb{N}}  B(y_n,{\rm dist}(A,B)+1/n).%
$$%
Since $\{ y_n\}$ is bounded we can assume it is weakly convergent to a
point $y$ in the closure of $B$. Then it must be the case that
$x,x^{\prime}\in B(y,{\rm dist(A,B)})$ from where, since $X$ is
strictly convex, $x=x^{\prime}$.

Therefore we can assume that $\{ x_n\}$ does not have
any convergent subsequence and so it is a separated sequence. Let $\varepsilon >0$ such that $d(x_n,x_m)\geq
\varepsilon$ for every $n\neq m$. Since this sequence is bounded and $X$ is
reflexive, we can also assume $\{ x_n\}$ is weakly convergent to a point $x$. Now we only have to apply
the UKK property in a similar way as the uniform convexity was applied in the
previous proposition to deduce that $x$ is in the closure of $A$ but
dist$(x,B)<{\rm dist}(A,B)$, contradicting the definition of dist$(A,B)$. $\Box$\vskip.5cm

The UKK property for $\Delta$-convergent sequences has been recently studied in \cite{esfe, kipa} for CAT$(k)$ spaces. If we assume that $X$ is a geodesic space such that bounded sequences have a unique asymptotic center which belongs to the convex hull of the sequence (see any of \cite{esfe, kipa} for definitions),
the previous proposition finds a metric counterpart that we state next and which can be proved exactly the same.

\begin{pro}
Let $X$ be a geodesic metric space with the UKK property
for $\Delta$-convergent sequences and the above-mentioned property for bounded sequences. Suppose also that the function $\mu (r,\varepsilon)$ given by the UKK property decreases with respect to the radius, then,
for $A,B\subseteq X$ nonempty and convex, it is the case that
$(A,B)$ has the property WUC.
\end{pro}

\begin{re}
Although the situation for Banach spaces is clear in the sense that uniform convexity
implies property UKK, the same seems far to be the case for geodesic spaces as defined for $\Delta$-convergent sequences. Actually
the UKK property for CAT$(k)$ spaces as shown in \cite{esfe, kipa} seems to be more connected with
the so-called Opial condition than with the uniform convexity. It is worth to recall at this point that
only Hilbert spaces and the spaces of sequences
$\ell_p$ are known to enjoy the Opial property. For more on this and related topics the interested reader
can consult Chapters 3, 4, 5 and 16 in \cite{kisi} or \cite[p. 102]{aydo}.
\end{re}

Next we show that WUC is enough to lead to a best proximity point
for a cyclic contraction. Due to notation purposes, we will denote $r$ as the contractive constant in the definition of cyclic contraction for the remainder of this section.

\begin{th}\label{teorema1}
Let $(X,d)$ be a metric space and $A$ and $B$ two nonempty subsets
of $X$ such that $(A,B)$ satisfies the property WUC. Assume that $A$
is complete. Let $T$ be a cyclic contraction on $A\cup B$. Then $T$
has a unique best proximity point $z$ in $A$ and the sequence $\{
T^{2n}x\}$ converges to $z$ for every $x\in A$.
\end{th}

{\it Proof.} As in \cite{suzu} we consider $d^*(x,y)=d(x,y)-{\rm
dist}(A,B)$. Then $d^*(Tx,Ty)\le rd^*(x,y)$ for $x\in A$ and $y\in
B$. In consequence $d^*(T^2x,Tx)\le rd^*(x,Tx)$ and $d^*(Ty,T^2y)\le r
d^*(Ty,y)$ for any $x\in A$ and $y\in B$.

Fix $x\in A$, $n\in N$ and let $m=n+k$ with $k\in \mathbb{N}$. Then
$$
d^*(T^{2m}x, T^{2n+1}x)\le r^{2n}d^*(T^{2k}x,Tx)
$$
$$
\le r^{2n}\sup\{ d(Tx,T^{2k}x):\; k\in {\mathbb{N}}\}=r^{2n}M(x).
$$

Proposition 3.3 in \cite{elve} guarantees that $M(x)$ is finite for each $x$. Hence, given $\varepsilon >0$ and taking $n$ such that
$r^{2n} M(x)<\varepsilon$ we have that
\begin{equation}\label{colas}
T^{2m}x\in B(T^{2n+1}x, {\rm
dist}(A,B)+\varepsilon)%
\end{equation}
for $m\geq n$ and so, by the property WUC, $\{T^{2n}x\}$ is
convergent. Now the proof follows the same patterns than the proof
of Theorem 3 in \cite{suzu}. Let $z\in A$ be the limit of
$\{T^{2n}x\}$, then
$$%
d^*(z,Tz)=\lim_{n\to \infty} d^*(T^{2n}x,Tz)\le \lim_{n\to\infty} r d^*(z, T^{2n-1}x)%
$$%
$$%
\le \lim_{n\to\infty} r(d(z,T^{2n}x)+d^*(T^{2n}x,T^{2n-1}x))%
$$%
$$%
\le \lim_{n\to\infty} r(d(z,T^{2n}x)+r^{2n-2}d^*(T^{2}x,Tx))=0.%
$$%
Thus, $z$ is a best proximity point. Now, we note that
$$
d(T^2z,Tz)\le d(z,Tz)={\rm dist}(A,B),
$$
and so $Tz$ is a best proximity point in $B$. Moreover, it follows
from the WUC that $T^2z=z$, and so $z$ is a fixed
points for $T^2$

Let $z^{\prime}$ be another proximity point in $A$, which will also be a fixed point for $T^2$. Then
$$%
d^*(z,Tz^{\prime})=d^*(T^2z,Tz^{\prime})%
$$%
$$%
\le r d^*(z^{\prime},Tz)=rd^*(T^2z^{\prime},Tz)\le
r^2d^*(z,Tz^{\prime}).%
$$%
Hence, $d^*(z,Tz^{\prime})=0$ and so $d(z,z^{\prime})=0$. $\Box$
\vskip.5cm

Still one further weakening of property WUC is possible to obtain a best proximity point result.

\begin{de}\label{wwuc}
Let $A$ and $B$ be nonempty subsets of a metric space $(X,d)$. Then
$(A,B)$ is said to satisfy the property W-WUC if for any $\{
x_n\}\subseteq A$ such that for every $\varepsilon >0$ there exists
$y\in B$ satisfying that $d(x_n,y)\le {\rm dist}(A,B)+\varepsilon $
for $n\geq n_0$ then there exists a convergent subsequence $\{ x_{n_k}\}$ of $\{ x_n\}$.
\end{de}

After this definition the following theorem is possible.

\begin{th}\label{322}
Under the same conditions of Theorem \ref{teorema1} with $(A,B)$ satisfying the property W-WUC assume that $A$ is a Chebysev set with respect to $B$,
it follows that $T$ has a unique best proximity point $z$ in $A$ and the sequence $\{
T^{2n}x\}$ converges to $z$ for every $x\in A$.
\end{th}

{\it Proof.} Following the same steps as in the previous proof we get that every subsequence of $\{ T^{2n}x\}$ has a convergent subsequence. Consider therefore
a convergent subsequence $\{ T^{2n_k}x\}$ of $\{ T^{2n}x\}$. Proceeding in a similar way, consider
$z\in A$ as the limit of $\{T^{2n_k}x\}$. Then
$$%
d^*(z,Tz)=\lim_{k\to \infty} d^*(T^{2n_k}x,Tz)\le \lim_{k\to\infty} r d^*(z, T^{2n_k-1}x)%
$$%
$$%
\le \lim_{n\to\infty} r(d(z,T^{2n_k}x)+d^*(T^{2n_k}x,T^{2n_k-1}x))%
$$%
$$%
\le \lim_{n\to\infty} r(d(z,T^{2n_k}x)+r^{2n_k-2}d^*(T^{2}x,Tx))=0.%
$$%
Thus, $z$ is a best proximity point. In the same way as above, but using the Chebyshev condition instead of property WUC, we obtain that $z$ is a fixed point of
$T^2$ and so uniqueness follows as in the previous proof. Also, since $T$ is relatively nonexpansive, we have that $\{ d(T^{2n}x,Tz)\}$ is a decreasing sequence with
$$
\lim_{k\to\infty}d(T^{2n_k}x,Tz)={\rm dist}(A,B).
$$
Therefore the tales of $\{ T^{2n}x\}$ are contained in $B(Tz, {\rm dist}(A,B)+\varepsilon)$ for $\varepsilon >0$.
Finally, the Chebyshev character of $A$ and the above equality implies that any convergent subsequence of $\{T^{2n}x\}$ must converge to $z$, which complete the proof of the theorem. $\Box$ \vskip.5cm

\begin{re}
The condition of $A$ to be Chebysev with respect to $B$ can be weakened to that one Chebyshev for proximinal points as introduced at the beginning of this section.
\end{re}





\begin{re}
The contractive condition imposed in the main result of \cite{suzu} is actually different than the one we have worked with in this work. More precisely, in \cite{suzu} it is supposed that there exists $k\in [0,1)$ such that
$$
d(Tx,Ty)\le k \max\{ d(x,y), d(x,Tx), d(Ty,y)\}\; +\; (1-k)\; {\rm dist}(A,B)
$$
for every $x\in A$ and $y\in B$. The proofs of our results in this section can be written, however, under this more general contractive condition without major modifications.
 \end{re}

\section{Reflexivity and convergence}\label{sec}
In this section we study the relation between cyclic contractions
and the semimetric $d_1$ defined in \cite{espi} for certain Banach
spaces. As a result we partially answer a question raised in \cite{elve} in the positive.



\begin{de} A pair $(A,B)$ of subsets of a metric space is said to be proximinal
if for each $(x,y) \in A\times B$ there exists $(x^{\prime},y^{\prime}) \in A \times B$ such that
$$d(x,y^{\prime})=d(x^{\prime},y)=\text{dist}(A,B).$$

If, additionally, we impose the condition that the pair of points
$(x^{\prime},y^{\prime}) \in A \times B$ is unique for each $(x,y) \in (A,B)$, then
we say that the pair $(A,B)$ is a sharp proximinal pair.
\end{de}

It was shown in \cite{espi} that when a pair of subsets $(A,B)$ of a strictly convex
Banach space is proximinal then the sets $A$ and $B$ also satisfy the following definition.

\begin{de} Let $A$ and $B$ be nonempty subsets of a Banach space $X$.
We say that $A$ and $B$ are proximinal parallel sets if the following
two conditions are fulfilled:
\begin{itemize}
\item [(1)] $(A,B)$ is a sharp proximinal pair.
\item [(2)] $B = A + h$ for a certain $h \in X$ such that $\|h\| = \text{dist}(A,B)$.
\end{itemize}
\end{de}

\begin{re}\label{trioparalelo}  The notation introduced in this remark will be used
along this section. Notice that, in the case of proximinal parallel sets, $a^{\prime}=a + h$ for $a \in A$ and $b^{\prime}=b-h$ for $b \in B$ are the proximinal points from $B$ and $A$, respectively, to $a$ and
$b$. Given $A$ and $B$ proximinal parallel sets in a Banach space $X$,
we will consider the set $C= A + 2h$, where $h \in X$ is such that
$B=A+h$. It is immediate to see that $A,B$ and $C$ are pairwise
proximinal parallel sets.
\end{re}
From now on, we will say that a pair $(A,B)$ satisfies a property if
each of the sets $A$ and $B$ has that property. We will need the
following technical result.

\begin{lem}\label{imagen} Let $A$ and $B$ be proximinal parallel subsets of a Banach space and
$T$ a cyclic contraction map defined on $A \cup B$. Then $T(a+h)=Ta - h$ for any $a \in A$ and $T(b-h)=Tb + h$
for any $b \in B$.
\end{lem}

{\it Proof.} Given $a\in A$, let $a^{\prime}=a+h$ be its proximinal point in
$B$. Since $T$ is relatively nonexpansive,
$$
\| Ta - Ta^{\prime}\| \leq \|a-a^{\prime}\|=\text{dist}(A,B),
$$
hence, by the uniqueness of the proximinal points, $T(a+h)=Ta-h$. $\Box$ \vskip.5cm

Next we define the semimetric $d_1$ introduced in \cite{espi}.

\begin{de}\label{d1}
Let $A,B$ and $C$ be as in Remark \ref{trioparalelo}. We consider
the function $d_1 : B \times B \rightarrow [0,\infty)$ by
$$%
d_1(x,y)=\inf \{r > 0 : y \in B(x-h,d+r) \cap B(x+h,d+r)\}.%
$$%
where $d=\text{dist}(A,B)$.
\end{de}

The following proposition follows in an immediate way.

\begin{pro}\label{cont1}
\begin{itemize}
\item [(1)] $d_1$ defines a semimetric on $B$.
\item [(2)] For every $x,y \in B$, $d_1(x,y) \leq \|x-y\|.$
\end{itemize}
\end{pro}

The next corollary immediately follows from (2) in the above
proposition.

\begin{co}
If $B(x,r)$ and $B_1(x,r)$ denote respectively the closed balls with
respect to the norm and the semimetric $d_1$ in $B$, then
$$%
B(x,r) \subseteq B_1(x,r)%
$$%
for any $x\in B$ and $r\geq 0$.
\end{co}

The next proposition gives sufficient conditions for the reverse
contention to happen in a nice way.

\begin{pro}\label{cont2}
Let $X$ be a reflexive and strictly convex Banach space which has
property $(H)$. Suppose that $A,B$ and $C$ are subsets of $X$ as
in Remark \ref{trioparalelo}. If $B$ is nonempty closed and
convex, then there exists $f:(0,\infty)\to (0,\infty)$ such that
$\lim_{r\to 0}f(r)=0$ and
$$%
B_1(x,r) \subseteq B(x,f(r)).%
$$%
\end{pro}

{\it Proof.} The existence of $f(r) >0$ such that $B_1(x,r)
\subseteq B(x,f(r))$ is immediate because $B_1(x,r)$ is bounded in
$X$. We need to prove that $f$ can be chosen so that the limit condition holds.
Since the balls $B_1$ are monotone with respect to the radius, it is enough to see that
$$%
\lim_{n\to \infty} f\left(\frac{1}{n}\right)= 0.%
$$%
Consider the set $B_1(x,\frac{1}{n})$ and let $x_n, y_n \in B$ such
that $d(x_n,y_n)=\text{diam}(B_1(x,\frac{1}{n}))$, where the diameter is taken with
respect to the norm metric (we assume the diameter is reached for simplicity).
Obviously, these points must belong to $\partial B_1(x,\frac{1}{n})$ where this
border is with respect to the topology induced in $B$ by the norm of $X$.
Since $\{B_1(x,\frac{1}{n})\}_n$ is a decreasing sequence of sets,
we have that $x_n$ and $y_n$ are in $B_1(x,1)$ for any $n$. Now, since
$B_1(x,1)$ is bounded closed and convex in a reflexive space $X$,
there exist subsequences $x_{n_k}$ and $y_{n_k}$ such that
$x_{n_k} \rightharpoonup z$ and $y_{n_k} \rightharpoonup w$, for
some $z, w \in  B_1(x,1)$. Moreover, since we know that $x_n \in
B_1(x,\frac{1}{n_0})$ for all $n \geq n_0$, we get that $z,w \in
B_1(x,\frac{1}{n})$ for all $n \in \mathbb{N}$. Hence it must be the case $z=w=x$, and so
the sequences $\{ x_n\}$ and $\{ y_n\}$ are weakly convergent
themselves to $x$.
Let $B^+=B(x+h,d)$, where $d=\| h\|={\rm dist}(A,B)$. Consider now the sequence $\{ P_{B^+}(x_n)\}$ of the radial
projections of $\{ x_n\}$ onto $B^+$. Thus,
$$
\|(x+h) - x_n\| =  \|(x+h) - P_{B^+}(x_n) \| + \|P_{B^+}(x_n) -
x_n\| = d + \|P_{B^+}(x_n) - x_n\|
$$
Since $x_n \in B_1(x,\frac{1}{n})$, it is the case that $\|(x+h) -
 x_n\| \leq d+ \frac{1}{n}$.
Hence, we obtain that $P_{B^+}(x_n) \rightharpoonup x$. Now,
property $(H)$ allows us to assure that this last sequence is also
convergent in norm, and, therefore, we also get that
$x_n \rightarrow x$. We can proceed analogously to proof that $y_n
\rightarrow x$. This implies that $d(x_n,
y_n)=\text{diam}(B_1(x,\frac{1}{n})) \rightarrow 0$, and so the
theorem follows by setting $\displaystyle f(1/n)=\text{diam}(B_1(x,\frac{1}{n}))$. $\Box$ \vskip.5cm

\begin{co}
Under the assumptions of the above proposition, the metric induced by the norm on $B$ and the semimetric $d_1$ are compatible in the sense of Definition \ref{compatible}.
\end{co}

Now we prove the main result of this section. The proof is inspired by those first appeared in
\cite{elvek, espi}.

\begin{th} \label{mainte}Let $A$ and $B$ be nonempty closed and convex subsets of
a reflexive and strictly convex Banach space $X$. Let $T : A \cup B
\rightarrow A \cup B$ be a cyclic contraction map. Suppose
that $X$ has the property $(H)$. 
Then there exists a unique point $b_0 \in B$ such that
$d(b_0,Tb_0)=\text{dist}(A,B).$ Moreover, there exists $h\in X$
and $B_0 \subseteq B$ such that, if $T^{\prime}(b)=Tb+h$ for $b \in B_0$,
then
\begin{itemize}
\item[(i)] $T^{\prime} : B_0 \rightarrow B_0$,
\item [(ii)] $b_0=Tb_0 + h$, and
\item [(iii)] $(T^{\prime})^n(b) \rightarrow b_0$ for each $b \in B_0$.
\end{itemize}
\end{th}
{\it Proof.} Given the pair $(A,B)$, let $A_0$ and $B_0$ be the subsets defined as follows:
$$A_0=\{x \in A : d(x,y^{\prime})= \text{dist}(A,B) \text{ for some } y^{\prime} \in B\},$$
$$B_0=\{y \in B : d(x^{\prime},y)= \text{dist}(A,B) \text{ for some } x^{\prime} \in A\}.$$

From the reflexivity of the space and the fact that $(A,B)$ is a closed and convex pair, the
pair $(A_0,B_0)$ is nonempty closed and convex itself. It follows from their definitions that the pair $(A_0,B_0)$ is proximinal with
$$
\text{dist}(A_0,B_0)=\text{dist}(A,B).
$$

Now we see that the mapping $T$ is still a {\it cyclic contraction}
on $A_0 \cup B_0$. Given $x_0 \in A_0
\subseteq A$, we have that there exists $y^{\prime}\in B_0$ such that
$d(x_0,y^{\prime})=\text{dist}(A,B)$. Since $T$ is a cyclic contraction on
$A \cup B$,
$$d(Tx_0,Ty^{\prime}) = \text{dist}(A,B).$$
Since $Ty^{\prime} \in A$, we get that $Tx_0 \in B_0$, and therefore $T(A_0)
\subseteq B_0$. In the same way we prove that $T(B_0) \subseteq
A_0$. Then $T : A_0 \cup B_0 \rightarrow A_0 \cup B_0$.
Moreover, since $\text{dist}(A_0,B_0)=\text{dist}(A,B)$,
$$
d(Tx,Ty) \leq kd(x,y)+(1-k){\rm dist}(A_0,B_0),
$$
for all $x \in A_0, y \in B_0$.

From Lemma $3.1$ of \cite{espi} we also have that $A_0$ and $B_0$ are proximinal
parallel sets. Let, therefore, $h \in X$ such that $B_0=A_0+ h$. This equality
of sets directly implies $(i)$. As it is shown in Lemma
\ref{imagen}, we have that for any $a\in A_0$ and $b \in B_0$,

$$T(a+h)=Ta - h \text{\espa and \espa} T(b-h)=Tb + h.$$

Now we prove that the mapping $T^{\prime}$ is a contraction on $B_0$ with respect to the semimetric $d_1$.
Let $x,y \in B_0$. Denote $r=d_1(x,y)$ and
$d=\text{dist}(A_0,B_0)$. To see that $T^{\prime}$ is a contraction with constant $k$ with respect to
$d_1$ it will be enough to show that
$$
T^{\prime}y \in B(T^{\prime}x - h, d + kr) \cap B(T^{\prime}x + h, d + kr).
$$
To see this, we proceed as follows:
$$T^{\prime}y \in B(T^{\prime}x - h, d + kr) \Leftrightarrow \|T^{\prime}y-(T^{\prime}x-h)\| \leq d+kr \Leftrightarrow$$
$$ \Leftrightarrow \|Ty+h-(Tx+h-h)\| \leq d+kr \Leftrightarrow  \|Tx-(Ty+h)\| \leq d+kr  \Leftrightarrow $$
\begin{equation}\label{siysolosi}
 \Leftrightarrow \|Tx-T(y-h)\| \leq d+kr.
\end{equation}
Since $T$ is a cyclic contraction,
$$ \|Tx-T(y-h)\| \leq k  \|x-(y-h)\| + (1-k)d \leq k(d+r)+(1-k)d=d+kr,$$
and hence the chain of inequalities (\ref{siysolosi}) holds.
To show that $T^{\prime}y \in B(T^{\prime}x + h, d + kr)$ we need to introduce a new mapping.
Let $\hat{T} :
B \cup C \rightarrow B \cup C$ the mapping defined by
$\hat{T}(b)=T(b)+2h$ if $b\in B$ and $\hat{T}(c)=T(c-2h)$ when
$c\in C$. It is immediate that $\hat{T}(B) \subseteq C$ and
$\hat{T}(C) \subseteq B$. Moreover, since
 $T(b)+h=T(b-h)$ and
$T(c-2h)=T(c-h)+h$, we have
$$
\|\hat{T}(b)-\hat{T}(c)\| =
\|T(b)+2h-T(c-2h)\| = \|T(b-h)-T(c-h)\|
$$
$$
\leq k \|b-c\|+(1-k)\; \text{dist}(A,B)= k \|b-c\|+(1-k)\; \text{dist}(B,C).
$$
Thus, $\hat{T}$ is a cyclic contraction. Now, repeating the same
reasoning as above we finally obtain that $T^{\prime}$ is a
$d_1$-contraction.

We conclude the proof of the theorem by applying Theorem \ref{convergencia}
to the set $B_0$ with the mapping $T^{\prime}$. The best proximity point $b_0$ is
given by $T^{\prime}$, so it is immediate that $(ii)$ and $(iii)$ hold.
The fact that $b_0$ is unique (as a best proximity point) even in
$B$ comes directly
by following similar patterns as in Theorem \ref{teorema1}. $\Box$
\vskip.5cm


We state the convergence of the iterates for any $b\in B$ in the following corollary.
\begin{co}\label{mainco}
Let $X$, $A$, $B$ and $T$ be as in the above theorem. Then, for every $b \in B$ the sequence of the iterates $\{T^{2n}(b)\}$ converges to $b_0$, the unique best proximity point of $T$ in $B$.
\end{co}
{\it Proof.} Let $b_0$ be the unique best proximity point of $T$ in $B$ given by the above theorem and consider $b\in B$. We proceed as follows:
$$
{\rm dist}(A,B) \leq \|T^{2n}(b)-T(b_0)\| =\| T^{2n}(b)-T^{2n-1}(b_0)\|
$$
(by iterating the contractive condition)
$$
\le k^{2n-1} \|T(b)-b_0\| + (1-k)\; {\rm dist}(A,B) \bigg{(}\sum_{i=0}^{2n-2} k^i\bigg{)} \leq  k^{2n-1} \|T(b)-b_0\| + {\rm dist}(A,B),
$$
therefore
\begin{equation}\label{ennorma}
\|T^{2n}(b)-T(b_0)\| \rightarrow {\rm dist}(A,B)
\end{equation}
as $n\to\infty$.

On the other hand, since $B$ is closed and convex in a reflexive space $X$ and $\{T^{2n}(b)\}$ is bounded, there exists a subsequence $\{T^{2n_{i}}(b)\}$ weakly converging to a point $x \in B$. Moreover, since
$$
\|x- T(b_0)\| \leq \liminf   \|T^{2n_{i_j}}(b)- T(b_0) \| ={\rm dist}(A,B),
$$
we get that $x=b_0$. Therefore, since $\{T^{n_i}b\}$ is any convergent subsequence, we obtain that
\begin{equation}\label{debil}
T^{2n}(b) \rightharpoonup b_0.
\end{equation}

Consider the closed ball $\hat{B}=B(T(b_0), {\rm dist}(A,B))$ and the radial projection of $T^{2n}(b)$ onto $\hat{B}$, which we denote by $P_{\hat{B}}(T^{2n}(b))$. Then
$$
\|T^{2n}(b)-T(b_0)\| =\|T(b_0)-P_{\hat{B}}(T^{2n}(b))\| + \|P_{\hat{B}}(T^{2n}(b))-T^{2n}(b)\|,
$$
together with (\ref{ennorma}), gives that $$\|P_{\hat{B}}(T^{2n}(b))-T^{2n}(b)\| \rightarrow 0.$$
Then, by (\ref{debil}), we obtain that $P_{\hat{B}}(T^{2n}(b)) \rightharpoonup b_0$. Now it is enough to apply property $(H)$, following the same reasoning as the one used in Proposition \ref{cont2}, to conclude
that ${T^{2n}(b)}$ converges to $b_0$. $\Box$\vskip.5cm

\begin{re}
If $b \in B_0$ in the above proof, the convergence is in fact immediate because $T^{2n}(b)=(T^{\prime})^{2n}(b)$.
\end{re}

\begin{re}
 In \cite{elve} they raised question of whether a best proximity point exists when $A$ and $B$ are nonempty closed and convex subsets of a reflexive Banach space. The conjunction of Theorem \ref{mainte} and Corollary \ref{mainco} partially answers this question as they also request strict convexity and property (H). Remember, however, that, as it was noticed in the Remark \ref{max}, the strict convexity assumption is a natural one in this kind of problems.
\end{re}

\subsection{Appendix: The CAT$\mathbf{(0)}$ case.}

One of the most important examples of uniformly convex uniquely geodesic metric spaces are the so-called CAT$(0)$ spaces. A CAT$(0)$ space is a space of nonpositive curvature in the sense of Gromov and they are characterized by having thinner triangles than the comparison ones in the $2$-dimensional Euclidean space. For a proper definition and further properties, the reader can consult Chapter II.1 in \cite{brha}.
CAT$(0)$ spaces can be viewed as a metric analog to the Hilbert spaces in the classical theory of nonlinear analysis. Properties studied in Chapter II.2 of \cite{brha} support very strongly this idea. It was studied in \cite{elvek} that relatively nonexpansive mappings with $T(A)\subseteq A$ and $T(B)\subseteq B$ are actually nonexpansive mappings with respect to the norm (i.e., $\|Tx-Ty\|\le \| x-y\|$) when the Banach space in their theorems happens to be a Hilbert space. It was also shown in \cite{espi} that, under similar assumptions to those of Theorem \ref{mainte} in this work, the metric induced by the norm coincides with the $d_1$-semimetric on the set $B_0$. The purpose of this appendix is to show that the same stands for CAT$(0)$ spaces.

Let, therefore, $X$, $A$, $B$ and $T$ be as in Theorem \ref{mainte} with no further assumption on $X$ but the fact of being a complete CAT$(0)$ space (also called a {\it Hadamard space}). Suppose dist$(A,B)>0$, otherwise there is nothing to prove, and define the sets $A_0$ and $B_0$ as before by
$$
A_0=\{x \in A : d(x,y^{\prime})= \text{dist}(A,B) \text{ for some } y^{\prime} \in B\},
$$
$$
B_0=\{y \in B : d(x^{\prime},y)= \text{dist}(A,B) \text{ for some } x^{\prime} \in A\}.
$$

The properties of the pair $(A_0,B_0)$ are summarized next.

\begin{pro}
The pair $(A_0,B_0)$ is a nonempty, closed and convex pair in $X$. Furthermore, each point $b\in B_0$ can be joined through a geodesic segment of length $d={\rm dist}(A,B)$ to its proximinal point $b-h$ in $A_0$ and viceversa.
\end{pro}

{\it Proof.} That they are closed follows in a straightforward way from their definition and the fact that $A$ and $B$ are both closed. The fact that $A_0$ and $B_0$ are nonempty also follows in a similar way to the linear case under the assumption of reflexivity, since it is a very well-known fact (see \cite{esfe, kipa}) that decreasing sequences of nonempty bounded closed and convex subsets of a CAT$(0)$ space have nonempty intersection. Finally, the convexity of the sets $A_0$ and $B_0$ follows from the convexity of the metric of CAT$(0)$ spaces (see Proposition 2.2 in Chapter II.2 of \cite{brha}). $\Box$ \vskip.5cm

The next thing we need to do is to define the semimetric $d_1$ on $B_0$. We will define it in such a way that a third set $C_0$ is not needed.

\begin{de} \label{d12}
We define the function $d_1 : B_0 \times B_0 \rightarrow [0,\infty)$ by
$$%
d_1(x,y)=\inf \{r > 0 : y \in B(x-h,d+r)\text{ and } y-h\in B(x,d+r)\},%
$$%
where $d=\text{dist}(A,B)$.
\end{de}

\begin{re}
Notice that Definition \ref{d1} and Definition \ref{d12} coincide in linear spaces.
\end{re}

\begin{th}
The semimetric $d_1$ coincides with the metric $d$ induced by $X$ on $B_0$.
\end{th}

{\it Proof.} This result follows as an easy application of {\it The Flat Quadrilateral Theorem} (\cite[p.181]{brha}). Indeed, consider the four point $x, y, x-h$ and $y-h$. Then $x$ (respectively, $y$) is the proximinal point of $x-h$ (rep., $y-h$) in $B_0$, and viceversa. In consequence, the angles $\angle_x(x-h,y)$, $\angle_y(x,y-h)$, $\angle_{y-h}(x-h,y)$ and  $\angle_{x-h}(y-h,x)$ are all greater than or equal to $\pi/2$. Therefore the Flat Quadrilateral Theorem implies that the convex hull of the points $x, y, x-h$ and $y-h$ is isometric to a rectangle in the $2$-dimensional Euclidean space. Now, by the Pythagorean theorem, it is immediate to deduce that
$$
B_1(x,\sqrt{d^2+r^2}-d)=B_0\cap B(x-h,\sqrt{d^2+r^2})= B_0\cap B(x,r),
$$
as we wanted to proof. $\Box$ \vskip.5cm

We close this appendix by observing that it is also possible to show  that the mapping $T^{\prime}(b)=Tb+h$ for $b\in B_0$ is actually a contraction. To see this we just need to proceed as in the proof of Theorem \ref{mainte} and recall, at the proper moment, that the convex hull of the points $x, y, x-h$ and $y-h$ is actually a rectangle.

\vfill

\noindent Departamento de An\'alisis Matem\'atico \\
Facultad de Matem\'aticas\\
Universidad de Sevilla\\
P.O.Box: 1160\\
41080-Sevilla

\noindent {\it emails:} espinola@us.es, aurorafl@us.es

\end{document}